\documentclass[letterpaper, 10 pt, conference]{ieeeconf} 

\usepackage[utf8]{inputenc}
\IEEEoverridecommandlockouts 
\usepackage{cite}
\usepackage{blindtext}

\usepackage{booktabs}
\usepackage{arydshln} 

\usepackage{amsmath} 
\usepackage{bbm}
\usepackage{bm}
\usepackage{amssymb}
\usepackage{mathtools}

\usepackage{amsthm}

\usepackage{comment}

\theoremstyle{definition}

\newtheorem{assumption}{Assumption}

\newtheorem{problem}{Problem}

\usepackage{hyperref}
\hypersetup{
    colorlinks=true,
    linkcolor=black,
    filecolor=black,      
    urlcolor=blue,
    citecolor=blue,
    }
\usepackage{tikz}

\usepackage{cite}

\usepackage{todonotes}

\title{\LARGE \bf
Energy-Optimal Multi-Agent Navigation as a Strategic-Form Game
}

\author{Logan E. Beaver, \emph{Member, IEEE} 
\thanks{L.E. Beaver is with the Department of Mechanical and Aerospace Engineering, Old Dominion University, Norfolk, VA USA (email: lbeaver@odu.edu).}%
}

\begin{document}
\maketitle
\thispagestyle{empty}

\begin{abstract}
This extended abstracts presents a method to generate energy-optimal trajectories for multi-agent systems as a strategic-form game. 
Using recent results in optimal control, we demonstrate that an energy-optimal trajectory can be generated in milliseconds if the sequence of constraint activations is known a priori.
Thus, rather than selecting an infinite-dimensional action from a function space, the agents select their actions from a finite number of constraints and determine the time that each becomes active.
Furthermore, the agents can exactly encode their trajectory in a set of real numbers, rather than communicating their control action as an infinite-dimensional function.
We demonstrate the performance of this algorithm in simulation and find an optimal trajectory in $45$ milliseconds on a tablet PC.
\end{abstract}

\section{Introduction}

Motion planning for multiple agents is, in general, a challenging problem.
One major concern is that agents trajectories are planned simultaneously, but constraints, such as collision avoidance, require the agents to exchange trajectory information.
A number of heuristic solutions have been developed to overcome this challenge.
Some, such as imposing an agent sequence or priority \cite{Turpin2013ConcurrentRobots}, impose a centralized ordering on the agents to ensure constraint satisfaction; however, there are many problems where imposing an external sequence is provably sub-optimal \cite{Ma2019SearchingFinding}.
Other approaches may be overly conservative; decentralized model predictive control either restricts how much the planned trajectory can deviate each iteration \cite{dunbar2006distributed} or imposes sequencing constraints \cite{qian2015decentralized}. 
Meanwhile, robust control \cite{cheng2020safe} assumes a worst-case scenario that may significantly degrade performance.

Alternatively, one could use game-theoretic approaches to generate trajectories for multiple agents simultaneously \cite{wang2019game}.
The agent actions are their control trajectories; thus, this is a differential game where each agent optimizes over an infinite function space \cite{Bryson1975AppliedControl}.
Each agent also has a payoff function, which is an explicit function of the actions taken by all agents.
This means that, in the decentralized setting, each agent must communicate a message containing their current state and control profile to every other agent.

In this work, we employ an original technique to generate optimal trajectories based on the sequence of constraint activations \cite{Beaver2021DifferentialTime,beaver2023graph,beaver2023lqocp}.
Namely, we parameterize the optimal trajectory using a finite number of variables and prove that this parameterization is optimal.
This implies that only communicating the sequence of active constraints, and the time at which the constraints become active, is sufficient for the agents to compress their optimal trajectories into finitely many real numbers.
This significantly reduces the communication and computational burdens on agents, which reduces the energy and time required for agents to plan their trajectories.

\section{Optimization Problem}

Consider a set of $N$ agents indexed by the set $\mathcal{A} = \{1,2,\dots,N\}$.
To simplify our exposition we omit the explicit dependence of state and control variables on $t$ where it does not lead to ambiguity.
Here we consider each agent $i\in\mathcal{A}$ to have double integrator dynamics in the 2D plane,
\begin{equation}
    \dot{\bm{p}}_i = \bm{v}_i, \quad \dot{\bm{v}}_i = \bm{u}_i,
\end{equation}
where $\bm{x}_i = [\bm{p}_i^{\intercal}, \bm{v}_i^{\intercal}]^{\intercal}\in\mathbb{R}^{4}$ is the state of agent $i$--consisting of the stacked position and velocity vectors--and $\bm{u}_i\in\mathbb{R}^2$ is the control action.
The agents seek to minimize the $L^2$ norm of their control input, i.e.,
\begin{equation}
    J = ||\bm{u}_i||^2,
\end{equation}
over some time horizon $[t_0, t_f]\in\mathbb{R}$.

The environment contains $K$ obstacles, which we index by $\mathcal{O} = \{1, 2, \dots, K\}$.
Each obstacle $k\in\mathcal{O}$ is circular with a position $\bm{p}_k$ and radius $r_k$.
To ensure safety, we circumscribe each agent $i\in\mathcal{A}$ in a ball of radius $r_i$ and impose pairwise minimum distance constraints, i.e., 
\begin{equation} \label{eq:constraints}
\begin{aligned}
    (r_i + r_j)^2 - ||\bm{p}_i - \bm{p}_j ||^2 &\leq 0, \quad \forall i,j\in\mathcal{A} \times \mathcal{A}\setminus\{i\}, \\
    (r_i + r_k)^2 - ||\bm{p}_i - \bm{p}_k ||^2 &\leq 0, \quad \forall i,k\in\mathcal{A}\times\mathcal{O}.   
\end{aligned}
\end{equation}
Note that we use the equivalent squared form of the constraint to guarantee smoothness of the derivatives.

Finally, combining each of these pieces yields the optimal control problem that each agent solves to determine its control trajectory.
\begin{problem} \label{prb:ocp}
    For each agent $i\in\mathcal{A}$,
    \begin{align*}
        \min_{\bm{u}_i(t)} & \int_{t_0}^{t_f} ||\bm{u}_i(t)||^2 dt \\
        \text{subject to:}& \\
        \dot{\bm{p}} &= \bm{v}, \quad \dot{\bm{v}} = \bm{u},\\
        (r_i + r_j)^2 -& ||\bm{p}_i - \bm{p}_j||^2 \leq 0 \quad \forall j\in\mathcal{A}\setminus\{i\} \\
        (r_i + r_k)^2 -& ||\bm{p}_i - \bm{p}_k||^2 \leq 0 \quad \forall k\in\mathcal{O},\\
        \text{given:}&\\
        &\bm{x}(t_0), \bm{x}(t_f).
    \end{align*}
\end{problem}
For systems with integrator dynamics, Pontryagin's minimization principle yields a differential equation that is necessary for optimality \cite{Beaver2021DifferentialTime},
\begin{equation} \label{eq:ode}
    \sum_{n=0}^{2}(-1)^n\frac{d^n}{dt^n}\Big(\frac{\partial J}{\partial \bm{p}^{(n)}} + \bm{\mu}^{\intercal}\frac{\partial \bm{g}}{\partial \bm{p}^{(n)}}\Big) = 0,
\end{equation}
where $\bm{\mu}^{\intercal}$ is a vector of time-varying influence functions that act like Lagrange multipliers to ensure constraint satisfaction, $\bm{g}$ is the vector of safety constraints, and the superscript $^{(n)}$ denotes the $n^{\text{th}}$ time derivative.
For a system with a quadratic cost, the differential equation is both sufficient and necessary for optimality in the unconstrained case \cite{ledzewicz2022pitfalls}.
Furthermore, we note that in energy-minimizing systems the agents tend to only instantaneously activate the constraints \cite{beaver2023graph,beaver2023lqocp}, i.e., outside of the pathological cases, we expect the agents to only touch the constraint boundary instantaneously; they do not slide along the surface of obstacles.
This motivates our main working assumption for this work,
\begin{assumption} \label{smp:main}
    The optimal trajectory of the agents is a sequence of unconstrained dynamical motion primitive connected by optimal junctions.
\end{assumption}
Note that we only impose Assumption \ref{smp:main} to simplify trajectory generation.
This assumption can be relaxed by using \eqref{eq:ode} to determine all possible constrained motion primitives.
Switching between motion primitives still only requires information about the connecting junctions, and thus our approach is still valid.

Under Assumption \ref{smp:main}, we note that in the unconstrained case, i.e., $\bm{\mu} = \bm{0}$, the dynamical motion primitive is
\begin{equation}
    \begin{aligned} \label{eq:unconstrained}
        \bm{u}(t) &= 6 \bm{c}_1 t + 2 \bm{c}_2 \\
        \bm{v}(t) &= 3 \bm{c}_1 t^2 + 2 \bm{c}_2 t + \bm{c}_3\\
        \bm{p}(t) &= \bm{c}_1 t^3 + \bm{c}_2 t^2 + \bm{c}_3 t + \bm{c}_4,
    \end{aligned}
\end{equation}
where $\bm{c}_1$--$\bm{c}_4$ are $8$ unknown constants of integration.
If the unconstrained trajectory is feasible, then the constants of integration are determined from the $8$ boundary conditions by solving a linear system of equations,
\begin{equation} \label{eq:linear-eqns}
\Bigg(
    \begin{bmatrix}
        t_0^3 & t_0^2 & t_0 & 1 \\
        3 t_0^2 & 2 t_0 & 1 & 0 \\
        t_f^3 & t_f^2 & t_f & 1 \\
        3 t_f^2 & 2 t_f & 1 & 0
    \end{bmatrix}
    \otimes
    I_{2\times2}
    \Bigg)
    \begin{bmatrix}
        \bm{c}_1 \\ \bm{c}_2 \\ \bm{c}_3 \\ \bm{c}_4
    \end{bmatrix} =
    \begin{bmatrix}
        \bm{p}(t_0) \\ \bm{v}(t_0) \\ \bm{p}(t_f) \\ \bm{v}(t_f)
    \end{bmatrix},
\end{equation}
where $\otimes$ is the Kronecker product, and $I_{2\times2}$ is the $2\times2$ identity matrix.

If the unconstrained trajectory is not optimal, then the system must encounter a junction at some unknown time $t_1$, i.e., a constraint must become active instantaneously.
This introduces $9$ additional unknowns--$8$ constants of integration for the unconstrained motion primitive that follows the junction, and $1$ for the unknown junction time $t_1$.
The corresponding equations are:
\begin{enumerate}
    \item Continuity in the state: $4$
    \item Continuity in control: $2$
    \item Constraint satisfaction: $1$
    \item $\big(\bm{p}_i(t_1) - \bm{p}_j(t_1)\big)\cdot\big(\bm{v}_i(t_1) - \bm{v}_j(t_1)\big) = 0$: $1$
    \item Continuity in $\dot{\bm{u}}_i(t_1)\cdot\bm{v}_i(t_1)$: $1$
\end{enumerate}
Here $i$ is the current agent and $j$ denotes the obstacle or agent $i$ is in contact with at $t_1$.

The $9$ equations listed above come from the optimality conditions that describe the jump in the Hamiltonian and the dual variables at the constraint junction \cite{Bryson1975AppliedControl,beaver2023graph,Beaver2021DifferentialTime}.
Note if we parameterize the junction by a vector $\bm{\phi}_1 = [\theta_1, t_1]^{\intercal}$, then
\begin{equation} \label{eq:rot-matrix}
    \bm{p}_i(t_1) = \bm{p}_j(t_1) + \begin{bmatrix}
        \cos(\theta) & -\sin(\theta) \\
        \sin(\theta) & \cos(\theta)
    \end{bmatrix},
\end{equation}
where $j$ again indexes the obstacle or agent that $i$ interacts with at the junction corresponding to $t_1$.

Note conditions 1--2 and \eqref{eq:rot-matrix} are linear in the trajectory coefficients with the same form as \eqref{eq:linear-eqns}.
Then, given a guess for the parameter $\bm{\phi}_1$, we calculate the entire trajectory of agent $i$ by solving a linear system of equations.
Finally, using conditions 4 and 5, we can perturb our guess $\bm{\phi}_1$ until all of the optimality conditions are satisfied.
This holds for arbitrarily many constraint activations, for both obstacle and agent collision avoidance.
In this way, we can construct the optimal trajectory by optimally selecting the optimal collection of parameters $\{\bm{\phi}_1, \bm{\phi}_2, \dots, \bm{\phi}_n\}$ that correspond to a sequence of $n$ constraints becoming active.

\section{Game-Theoretic Formulation}

To solve our multi-agent path planning problem as a strategic-form game, we define the tuple \cite{Chremos2020SocialDilemma}
\begin{equation} \label{eq:game}
    \mathcal{G} \coloneqq (\mathcal{A}, \mathcal{M}, J),
\end{equation}
where $\mathcal{A}$ is the set of agents, $\mathcal{M} = \{\bm{\Phi}_i, i\in\mathcal{A}\}$ is a \emph{message} that encodes agent's trajectory, and $\bm{J} = \{J_i, i\in\mathcal{A}\}$ is a set of payoff functions for each agent.
Each message $M_i$ contains the sequence of constraint activations and the corresponding junctions $\phi_k, k = 1,2,\dots,n$.
The payoff function $J_i : \mathcal{M}\to\mathbb{R}\cup\{\infty\}$ takes all agent actions as inputs and yields the cost of agent $i$.
Note that $J_i$ is exactly equal to the solution of Problem \ref{prb:ocp} if the resulting trajectory is feasible and infinity otherwise.
This defines our strategic-form game, where the Nash equilibrium is a collection of agent trajectories that solve Problem \ref{prb:ocp}.
However, the message $M_i$ of each agent effectively compresses the entire optimal trajectory down to finitely many real numbers.

While the solution using Pontryagin's minimum principle is, in general, only necessary and sufficient to find critical points in the optimization space \cite{ledzewicz2022pitfalls}, the only critical point is the global maximum in problems with linear dynamics, a quadratic cost, and no constraint activations.
This insight, along with our proposed approach game-theoretic formalization, opens the door to three significant research problems:
\begin{enumerate}
    \item Is Pontryagin's minimum principle sufficient for global optimality with quadratic constraints of the form \eqref{eq:constraints}?
    \item There are many feasible sequences of constraint activaitions, are there efficient methods to find the global sequence of constraint activations?
    \item What conditions are sufficient for a Nash equilibrium to exist for the strategic form game \eqref{eq:game}?
    \item How can the cost function $J_i$ be designed, i.e., using \textit{mechanism design} to ensure that there is a unique dominant strategy equilibrium \cite{Dave2020SocialMedia}?
\end{enumerate}
Rather than addressing the above research directions on solution existence, we demonstrate that feasible solutions can be found efficiently for multi-agent systems in obstacle-filled environments.

\section{Numerical Results}

\begin{figure}[ht]
    \centering
    \includegraphics[width=0.75\linewidth]{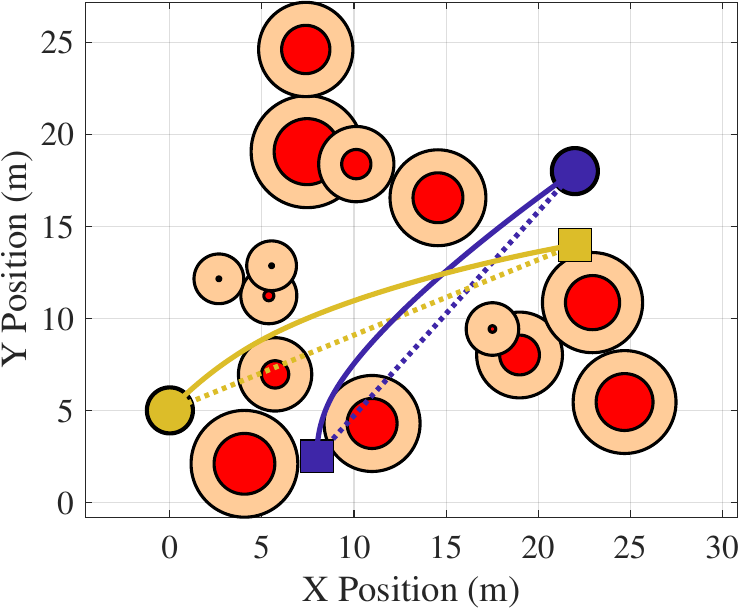}
    \caption{The two agents, blue and yellow, navigating from their initial state (circles) to goal states (squares) while avoiding obstacles. The red obstacles are inflated to compensate for the size of the agents during motion planning.}
    \label{fig:trajectories}
\end{figure}

To demonstrate our ability to rapidly generate trajectories, consider the randomly generated sphere world in Fig. \ref{fig:trajectories}.
It depicts two agents, with their initial positions marked with a circle and their final positions marked with a square.
Both agents start from rest; they must navigate through the environment while avoiding collisions with the spherical objects and each other.
Each agent is circumscribed in a circle of radius $r=1.25$ meters, and the cream colored halo around the spherical obstacles shows how they are inflated to compensate for the agents' size during motion planning.

First we solve for the unconstrained trajectory of each agent, denoted by the dashed line in Fig. \ref{fig:trajectories}.
Using \eqref{eq:unconstrained} with the boundary conditions takes approximately $1$ ms to generate on a tablet PC (i5-7300 @ 2.6 GHz, 8 GB of RAM).
Both unconstrained trajectories travel through an object, thus they are infeasible.
For the next step, each agent must construct a trajectory made of two unconstrained motion primitives.

\begin{figure}[ht]
    \centering
    \includegraphics[width=0.75\linewidth]{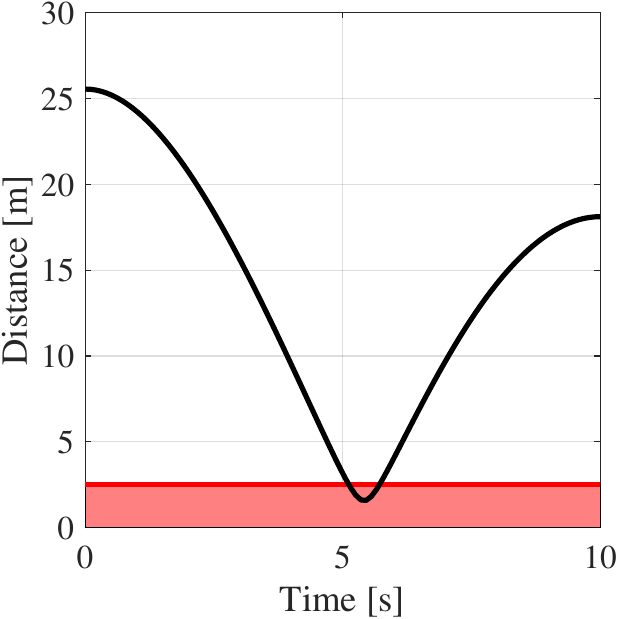}
    \caption{The pairwise safe distance constraint between the agents for the environment shown in Fig. \ref{fig:trajectories}.}
    \label{fig:safety}
\end{figure}

We parameterize the trajectory by $\sigma_i = [t_1, \theta_1]$, where $i\in\{1,2\}$ denotes the agent, $t_1$ is the time the agent arrives at the obstacle, and $\theta_1$ describes the point of contact between the agent and the obstacle.
We generate the optimal trajectories using Matlab's \textit{lsqnonlin} to minimize the least squares objective function,
\begin{equation}
\begin{aligned}
    J =& \Bigg(||\dot{\bm{u}}_i(t_1)\cdot\bm{v}_i(t_1)||\Bigg)^2 \\
    &+ \Bigg(\bm{v}_i(t_1)\cdot\begin{bmatrix}
        \cos(\theta_1) & \sin(\theta_1) \\
        -\sin(\theta_1) & \cos(\theta_1)
    \end{bmatrix}
    \begin{bmatrix}
        1 \\ 0
    \end{bmatrix}
    \Bigg)^2.
    \end{aligned} \label{eq:least-squares}
\end{equation}
We minimize \eqref{eq:least-squares} to a value of $10^{-7}$ in approximately $45$ milliseconds, which yields the optimal constrained trajectory for each individual agent--without considering the inter-agent collisions.
This is clear from Fig. \ref{fig:safety}, which shows the distance between the agents at each time.
The agents reach a maximum overlap of $0.4$ m at approximately $5.3$ seconds, which implies a collision between the agents.

\begin{figure}
    \centering
    \includegraphics[width=0.75\linewidth]{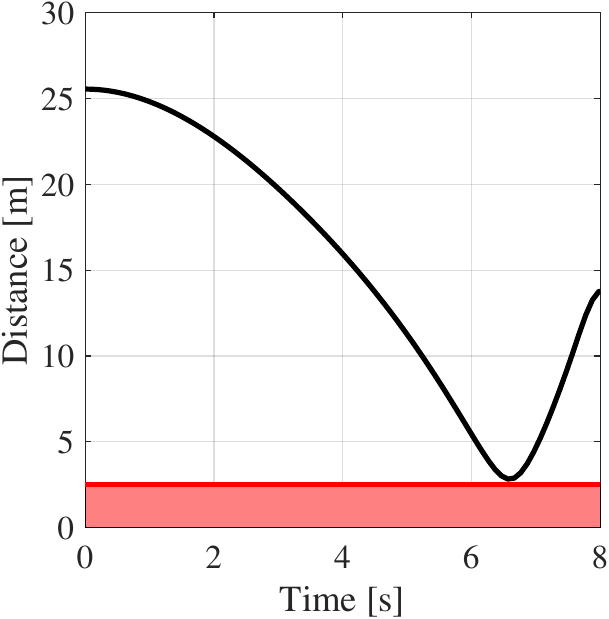}
    \caption{The relative distance between the agents satisfies the safety constraints when their arrival times are delayed or advanced to avoid collisions.}
    \label{fig:safety-delayed}
\end{figure}

Finally, there are several approaches to generate feasible agent trajectories that coordinate to avoid collisions:
\begin{enumerate}
    \item Change the final arrival time of each agent to ensure they cross at different times.
    \item Impose a delayed crossing time at the conflict point to ensure safety.
    \item Solve the centralized problem to guarantee safety.
\end{enumerate}
We opted for the first option, as it only takes $45$ ms to generate optimal constrained trajectories.
In a game-theoretic sense, we have added an additional parameter $t_i^f$ to each of the agents, that the agents must keep as closely as possible to the original time $t^f = 10$ sec while avoiding collisions.
The resulting arrival times are $t_1^f = 8, t_2^f = 12$ seconds for the blue and yellow agents, respectively.
The evolution of the pairwise distance between agents is shown in Fig. \ref{fig:safety-delayed}.

\section{Conclusions}
This extended abstract demonstrated how recent developments in optimal control can be used to generate trajectories for multiple agents.
First, we formulated the decentralized optimal control problem and parameterized the optimal trajectory.
We proved that our parameterization is optimal, and that multi-agent motion planning can be cast as a differential game.
Using our parameterization, the differential game is converted to a strategic-form game, where agents must select their optimal parameterizations to generate their trajectories.
We completed our analysis by demonstrating that safe optimal trajectories can be generated in milliseconds on a standard tablet PC.
Finally, we suggested several compelling areas of future research throughout this extended abstract.

\bibliography{mendeley,refs,IDS_Pubs}

\end{document}